\documentclass[12pt]{article}
\usepackage{amsfonts}
\usepackage{epsfig}
\title{Notes on Shapes of Polyhedra}
\author{Rich Schwartz}

\newtheorem{theorem}{Theorem}[section]

\newtheorem{lemma}[theorem]{Lemma}

\newtheorem{corollary}[theorem]{Corollary}

\def\startproof{{\bf {\medskip}{\noindent}Proof: }}

\def\endproof{$\spadesuit$  \newline}

\def\C{\mbox{\boldmath{$C$}}}%
\def\H{\mbox{\boldmath{$H$}}}%
\def\P{\mbox{\boldmath{$P$}}}%
\def\R{\mbox{\boldmath{$R$}}}%
\def\Z{\mbox{\boldmath{$Z$}}}%

\begin{document}
\maketitle

\section{Introduction}

Bill Thurston wrote a beautiful paper called {\it Shapes of
  Polyhedra\/}. I once lectured on this paper during a graduate class
I taught at the University of Chicago, and recently (Fall 2013) I
tried again during my graduate class at ICERM/Brown.  I found the
paper hard-going both times. In the intervening years, Thurston
published an updated and improved version, but I found the new version
hard going as well.

I wrote these notes for the ICERM class, and
some people (both in and out of class) found them very useful.
After some encouragement, I decided to
put them on the arXiv, so that they have
a public and stable home.
I am pretty sure that the proofs are
correct, but perhaps I am still missing
something.  Take them or leave them.

One of the most difficult parts of the paper
is the discussion of complex hyperbolic cone
manifolds.  For one thing, the definition is
hard to grasp.  For another thing, it is 
hard to see that the moduli spaces in question
really are  complex hyperbolic cone
manifolds according to the definition.
In these notes, I will explain things without
relying on cone manifolds at all.  Rather, I will introduce
related objects which are easier to understand.

At the end of this note, I'll list some other
references, one from Curt McMullen and several
from John Parker, which treat topics closely
related to Thurston's paper.

\subsection{Main Results}

Part of Thurston's paper deals with triangulations of 
the sphere and the other part deals with moduli spaces
of flat cone spheres. I'll talk about
the flat cone spheres first.  A {\it flat cone sphere\/} is a
metric on the sphere which is locally isometric to
the Euclidean plane except at finitely many points,
where it has positive conical singularities. This means
that a neighborhood of the point is isometric to a
Euclidean cylinder with cone angle $2\pi-\theta$
for some $\theta \in (0,2\pi)$.  The number
$\theta$ is called the cone deficit.  One might
say that the ordinary points have cone deficit zero.

Let $\theta_1,...,\theta_m$ be a finite list of
positive numbers such that $\sum \theta_i=4 \pi$.
Let ${\cal M\/}$ denote the moduli space of
similarity classes of
flat cone spheres with labeled cone deficits
$\theta_1,...,\theta_m$.  Sometimes we shall
take the ``labeled moduli space'', in which
all the cone points are labeled. In this
moduli space, two flat cone structures
are close if there is a near-isometry which
maps cone points to cone points and
respects the labels.  At other times, we shall
take the ``unlabeled moduli space''.
In this space, two flat cone structures are
close if there is a near isometry between
them which maps cone points to cone points
and respects the deficit values.
If all the cone deficits are distinct,
the two spaces are the same.  In general,
the unlabeled space is a quotient
of the labeled space by a finite group of
isometries.

Let $\C\H^n$ denote
complex hyperbolic space.
The first main result in Thurston's paper is

\begin{theorem}
\label{basicX}
Let $\cal M$ be the labeled moduli space.
$\cal M$ has a natural metric with respect
to which it is locally isometric to
$\C\H^{m-3}$.
\end{theorem}

\noindent
{\bf Remark:\/}
Theorem \ref{basicX} is not quite true in the
unlabeled case.  For instance, suppose
$m=4$ and $\theta_i=\pi$ for all $i$.
In this case, ${\cal M\/}$ is isometric
to the familiar modular orbifold.
\newline

The metric on $\cal M$ is incomplete whenever
there are two deficits $\theta_i$ and $\theta_j$
such that $\theta_i+\theta_j<2 \pi$. Really,
the example mentioned in the remark is the
only nontrivial example where the metric
is complete.

In the incomplete case, which essentially
always happens, 
Thurston goes on to prove (in some sense) that
the completion is a complex hyperbolic cone
manifold.  
The most interesting case occurs when the
list of deficits satisfies two additional
conditions.

\begin{enumerate}
\item If $\theta_i+\theta_j<2 \pi$ then
$2\pi$ is an integer multiple of $2\pi-\theta_i-\theta_j$.
\item If $\theta_i=\theta_j<\pi$ then
$2 \pi$ is an integer multiple of $\pi-\theta_i$.
\end{enumerate}
This case occurs for the flat cone spheres which
arise in connection with the triangulations.
In this case, Thurston proves a stronger result,
one highlight of the paper.

\begin{theorem}
\label{lattice}
Let $\cal M$ be the unlabeled moduli space.
If the deficit list satisfies the additional
conditions, then
there is a lattice $\Gamma$ acting on
$\C\H^{m-3}$ so that the metric completion
of $\cal M$ is isometric to $\C\H^{m-3}/\Gamma$.
\end{theorem}

\noindent
{\bf Remark:\/}
Theorem \ref{lattice} is not quite true
for the labeled space.  However,
Theorem \ref{lattice} is true for the labeled space
if all the cone
deficits are distinct, or if it never happens
that there are two equal cone deficits
less than $\pi$. 
\newline

Let ${\bf Eis\/}$ denote the Eisenstein lattice,
$\Z[\omega]$, where $\omega=\exp(2 \pi i/3)$ is the
usual cube root of unity.  The points in
${\bf Eis\/}$ are the vertices of the usual triangulation
of the plane by equilateral triangles.

Say that a triangulation of the sphere is {\it combinatorially positive\/}
if there are never more than $6$ triangles around a vertex.
Each combinatorial positive triangulation
gives rise to a flat cone sphere -- one just glues together
equilateral triangles in the same pattern.
If $k$ triangles go around a vertex, the corresponding
deficit is $(6-k)\pi/3$.  

Let's think of these triangulations as giving points
in the unlabeled moduli space ${\cal M\/}$. 
We call ${\cal M\/}$
{\it special\/} if it contains at least one point corresponding
to a triangulation.  Call a point in $\cal M$ a
{\it triangulation point\/} if it comes from a triangulation.
The following result is not explicitly stated in
Thurston's paper, but it is implied by other results.

\begin{theorem}
\label{tri1}
If ${\cal M\/}$ is special, the set of
triangulation points is dense in $\cal M$.
There is a single lattice
$\Gamma$, acting on $\C^{1,9}$, defined over
${\bf Eis\/}$, such that every special
moduli space is isometric to some
stratum of $\C\H^9/\Gamma$. 
\end{theorem}

The set ${\bf Eis\/}^{1,9}$ denotes the set of
vectors in $\C^{1,9}$ having coordinates in
${\bf Eis\/}$. 
Here is Thurston's main result about 
combinatorially positive triangulations,
the other highlight of the paper:

\begin{theorem}
\label{tri2}
There is a natural bijection between the set of
combinatorially positive triangulations and
the set of vectors in
${\bf Eis\/}^{1,9}/\Gamma$ having
positive square norm with respect to a
$\Gamma$-invariant Hermitian form
$H$ of type $(1,9)$.
Here $\Gamma$ is a lattice defined
over ${\bf Eis\/}$ and $H$ is also
defined over ${\bf Eis\/}$.
The square norm
$H(V,V)$ of a positive vector $V \in {\bf Eis\/}^{1,9}$
is $3$ times the number of
triangles in the triangulation corresponding to $V$.
\end{theorem}

\subsection{Organization of the Notes}

I'll explain the proof of Theorem \ref{basicX} in
\S \ref{easy}. The proof essentially follows
Thurston's outline, except that I do some of
the details differently.

Thurston proves Theorem \ref{lattice} in three steps.
\begin{enumerate}
\item Ignoring Conditions 2 and 3 above, 
the completion of $\cal M$ is always a finite volume
complex hyperbolic cone manifold.
Conditions 2 and 3 together imply that the
completion of $\cal M$
has codimension $2$ orbifold singularities.
\item A complex hyperbolic cone manifold with
codimension 2 orbifold singularities is in fact an
orbifold.
\item Every finite complex hyperbolic orbifold is
a lattice quotient.
\end{enumerate}

I'll try to explain Theorem \ref{lattice} in a
different (but related) way
which avoids the discussion of cone manifolds
and most of the discussion of orbifolds.
The objects I'll work with sound more
technical, but in fact they are easier
to understand because they refer very little
to the structure of the singular locus.
My route to Theorem \ref{lattice} 
works like this:

\begin{enumerate}
\item Ignoring Conditions 2 and 3 above, 
the completion of 
$\cal M$ is always a finite volume
complex hyperbolic stratified manifold
with a fibered cone structure.
Conditions 2 and 3 together imply that 
the completion of $\cal M$
has codimension $2$ orbifold singularities.
\item Theorem
\ref{strat} below: A complex hyperbolic stratified
manifold with a fibered cone structure and
codimension $2$ orbifold singularities
is a lattice quotient.
\end{enumerate}

In \S \ref{defs}, I'll explain the terms used above.
In \S \ref{struct}, I'll explain why one gets such
aobjects from the details in Thurston's paper.
In \S \ref{latq}, I'll prove Theorem \ref{strat},
thereby completing the proof of Theorem \ref{lattice}.

In \S \ref{triproof1} and \ref{triproof2}, I'll
prove Theorems \ref{tri1} and \ref{tri2} respectively.
These results are essentially interpretations
of Theorem \ref{lattice}.  We just have to look
back over the various constructions and see that
they give the statements in
Theorems \ref{tri1} and \ref{tri2}.

\section{Proof of Theorem \ref{basicX}}
\label{easy}

Let ${\cal M\/}$ be the labeled moduli space.

\subsection{Step 1: Spanning Trees and Triangulations}

A {\it embedded spanning tree\/} on a flat cone sphere
$\Sigma$ is spanning tree which has the cone points
as vertices and no crossing edges.  For instance,
on the regular tetrahedron, the three edges emanating
from a single vertex would be an embedded spanning tree.

\begin{lemma}
$\Sigma$ has an embedded spanning tree.
\end{lemma}

\startproof
First of all, $\Sigma$ does have a spanning tree.
One can connect every cone point to every other
one by some straight line segment, and now one
can choose a subgraph which is a spanning tree.
Let $\tau$ be the spanning tree of minimum
length.  If a pair of edges in $\tau$ cross,
then we can find a finite cycle
$e_1,...,e_k$ such that $e_1$ and $e_k$ cross.  
But then we can switch the crossing, as shown
in Figure 1.  The switch gives a shorter
spanning tree.
\endproof

\begin{center}
\resizebox{!}{1.5in}{\includegraphics{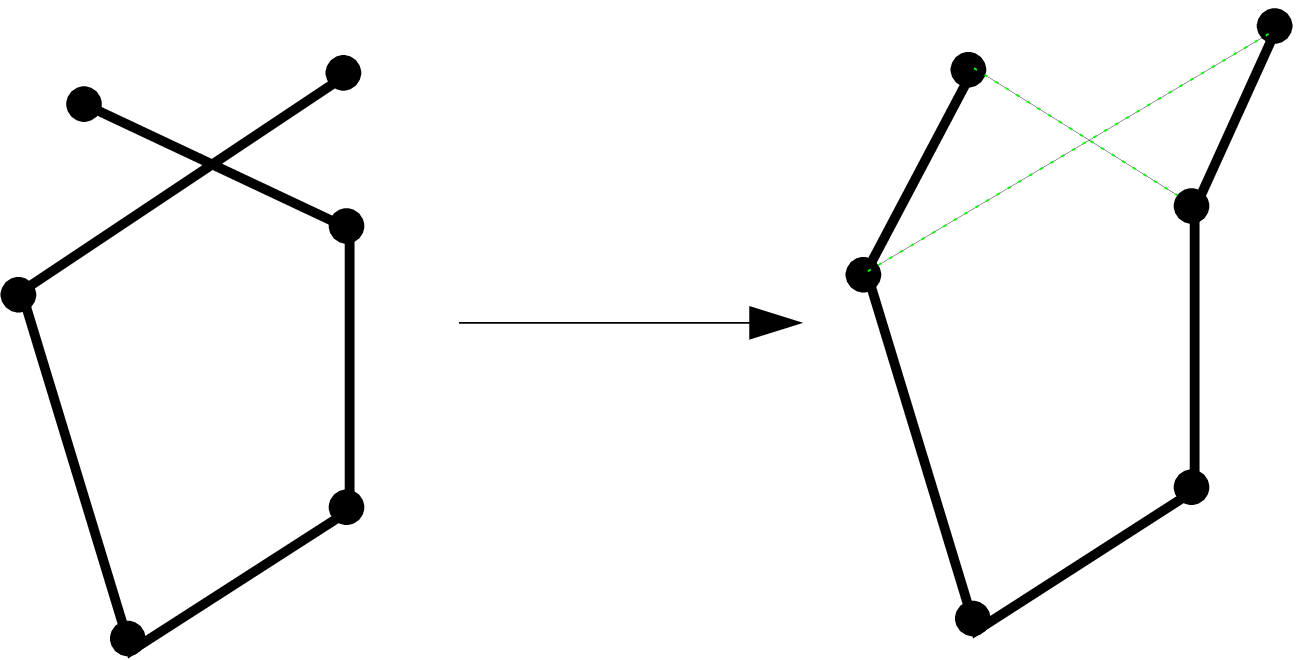}}
\newline
Figure 1: shortening the spanning tree.
\newline
\end{center}

For the proof of the next result, and for later
purposes, say that a {\it pseudo-polygon\/} is a flat metric on
a disk whose boundary is locally isometric to the
boundary of a polygon.  One typically gets a pseudo-polygon
by immersing a polygon in the plane and pulling back the
metric.

\begin{lemma}
A pseudo-polygon $P$ has a triangulation whose edges are
straight line segments and whose vertices are the vertices of $P$.
\end{lemma}

\startproof
The proof goes by induction on the number of edges of $P$.
The result is obvious if $P$ has $3$ edges.
If $P$ has more than $3$ edges, then let $v$ be a vertex of
$P$ and let $e$ be an incident edge.  Let $L_t$ be a family
of rays emanating from $v$ so that $L_0$ extends $v$ and
the initial portion of $L_t$ lies in $P$ for $t>0$ small.
Let $S_t \subset L_t$ be the longest initial portion of
$L_t$ contained in the interior of
$P$.  Note that $S_t$ is positive for
all $t \in (0,\theta)$ where $\theta$ is the interior
angle at $v$.  There must be some $s \in (0,\theta)$ such
that the endpoint of $S_s$ is another vertex.  But then
$S_s$ is an embedded segment connecting two distinct
vertices of $P$.  This segment divides $P$ into two
pseudo-polygons which both have triangulations by
induction.
\endproof

\begin{corollary}
\label{triangulate}
$\Sigma$ has a triangulation which just has the
cone points as vertices.
\end{corollary}

\startproof
Let $\tau$ be an embedded spanning tree in $\Sigma$.
The complement $\Sigma-\tau$ has a flat metric.
The completion of this metric is a pseudo-polygon
having twice as many edges as $\tau$.  But then
we can triangulate this pseudo-polygon.
\endproof

\noindent
{\bf Remark:\/}
Thurston proves Corollary \ref{triangulate} with
a canonical construction using Voronoi cells and
the dual Delaunay triangulation.  However, I 
found it hard to see why this construction
gives a triangulation rather than a union of
triangles, with pairwise disjoint interiors,
which perhaps only covers part of the surface.
This is why I prefer the spanning tree approach.
In class Saul Schliemer suggested that one
could remove the vertices, pass to the
universal cover, and then take the
Delaunay triangulation there.  That seems
to work more convincingly than the argument
in the paper, though I still
prefer the spanning tree approach.

\subsection{Step 2: Local Coordinates}
\label{loc}

Let $\Sigma$ be a flat cone sphere, a point in
$\cal M$.  We are really interested in flat
cone spheres mod similarity, but first we consider
the larger space of flat cone structures.
Let $\tau$ be an embedded spanning tree on $\Sigma$.
A small neighborhood in $\cal M$ consists of
flat cone spheres having an embedded spanning
tree combinatorially identical to, and nearby, $\tau$.

We orient the edges of $\tau$ in some way.
Let $P$ be the pseudo-polygon which is the
completion of $\Sigma-\tau$. The developing
map $D: P \to \C$ is well-defined because
$P$ is simply connected and has a flat metric.
We label each edge $e$ of $P$ by the complex number
\begin{equation}
f(e)=D(e_+)-D(e_-).
\end{equation}
Here $e_+$ is the head vertex of $e$ and
$e_-$ is the tail vertex.  Call this label $f(e)$.
If $e$ and $e'$ are the two edges glued together,
then we have a relation of the form
$f(e')=u_e f(e)$,
where $u_e$ is some unit complex number that only
depends on the list of cone deficits.  One computes
$u_e$ by taking a loop in $\Sigma$ which starts
and ends at (say) the midpoint of $e$ and avoids
$\tau$.  This loop encloses some number of cone
points, and the number $u_e$ is $\exp(i\theta_e)$
where $\theta_e$ is either the sum of the
cone deficits enclosed by the loop or $2\pi$ minus
that sum.  Which option depends on the orientation
of the loop.

These labels make sense on all flat cone
structures near $\Sigma$.
If we multiply all labels by some complex number
$\lambda$ we get the same structure up to similarity.
Moreover, the 
labels specify a pseudo-polygon which we can
then glue together to get a point in $\cal M$.
Thus, two nearby flat cone spheres are similar
to each other if and only if their labels
differ by this kind of scaling.  In short,
we have given local coordinate charts into
projective space $\C\P^{n-3}$.  Here
$n$ is the number of cone points.

Our labels extend to give a system of labels on 
a triangulation of $\Sigma$ extending $\tau$.
The labels on the remaining (oriented) edges
are linear combinations of the labels of the
edges of $\tau$.  

Suppose that we choose a different spanning tree.
Each edge in the new tree cuts through a finite
number of triangles of the old triangulation.
When we develop these triangles out into the
plane, we express the new edge label as some
complex linear combination of old labels.  Hence
the changes of coordinates are complex linear.
Remembering that we need to mod out by scaling,
we see that the overlap functions for our
charts are complex projective transformations.
In particular, $\cal M$ is a complex projective
manifold.

\subsection{Step 3: The Hermitian Form}

We fix a spanning tree $\tau$ on $\Sigma$ and
consider the local coordinates on the
edges of the pseudo-polygon $P$.  Call this
larger space $\cal P$.  A neighborhood
of $\cal M$ about $\Sigma$ is the quotient
of $\cal P$ by scaling, as discussed above.

We have the area function $A: {\cal P\/} \to \R$.
The area of a triangle spanned by edges $z$ and $w$ is
\begin{equation}
\label{her}
\pm \frac{i}{4}(z\overline w-w\overline z)
\end{equation}
The sign depends on whether the vectors $\{z,w\}$
make a positively oriented or a negatively oriented
basis.

When we express $A$ as a function of the
edge labels, we get a finite number of sums
of terms like the one in Equation \ref{her},
where $z$ and $w$ are various complex linear
combinations of the edge labels. Hence,
$A$ is the diagonal part of a Hermitian form.
The coordinate changes are isometries relative
to this form because changing the spanning
tree does nothing to the area.

Now I'll explain why the 
Hermitian form has type $(1,n-2)$.  But then
the space $\cal M$ locally has the structure
of $\C\P^{n-3}$.

To explain the type of the Hermitian form, suppose
that there are $2$ cone deficits, say 
$\theta_1$ and $\theta_2$ such that
$\theta_1+\theta_2<2\pi$.  Then we join
the corresponding cone points by a straight line
segment and slit $\Sigma$ open along this line segment.
We then glue in an appropriate portion of a cylinder
to produce a new flat cone sphere $\Sigma_{12}$ with one fewer
cone point.  Figure 2 shows a schematic view of this.

\begin{center}
\resizebox{!}{1.4in}{\includegraphics{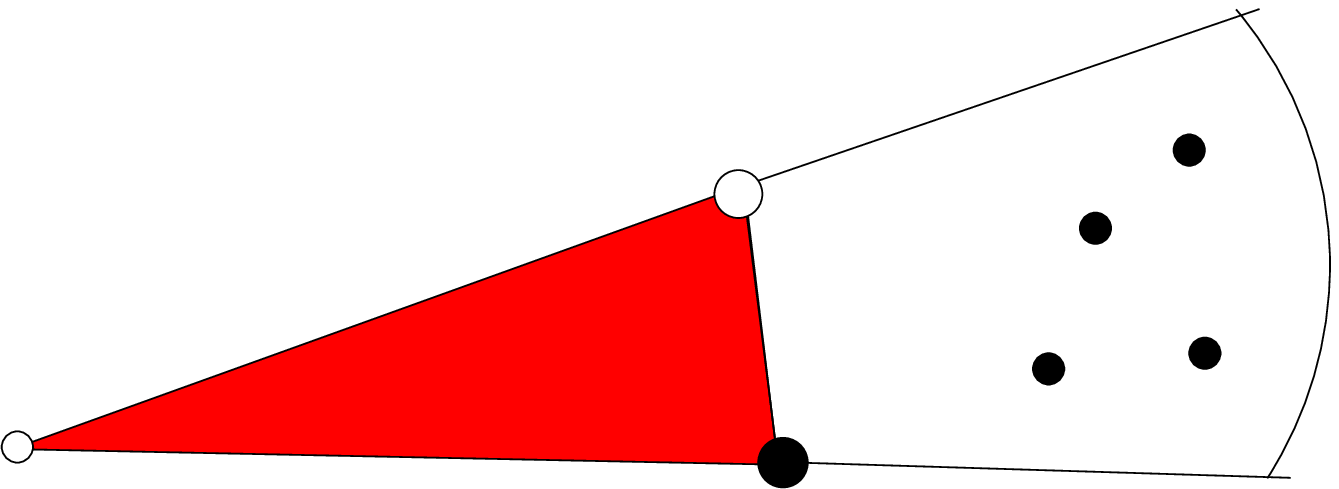}}
\newline
Figure 2: adding in the cone
\newline
\end{center}

Imagine that we have given linear coordinates
$w_1,...,w_{n-3}$ 
on $\Sigma'$.  These coordinates tell us how to
develop $\Sigma'$ out into the plane.  Let's
say that the apex of the added (red) cone goes
to the origin.  Then there is some complex
number $z$ which describes the position of one
image of $\theta_1$ under the developing map.
Fgure 3 shows what we are talking about.

\begin{center}
\resizebox{!}{2.6in}{\includegraphics{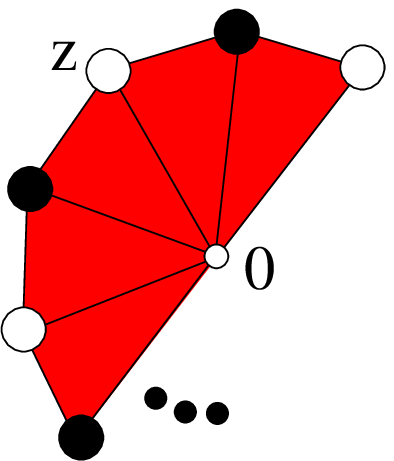}}
\newline
Figure 3: developing out the cone
\newline
\end{center}

Of course, the point $z$ depends on which
image under the developing map we choose. In
general, we have countably many choices.
However, once we make one choice for $\Sigma$,
we can make the same choice, so to speak, for
nearby structures.  The coordinate $z$ is
a complex linear function of the linear
coordinates on the space $\cal X$ described
above.  Thus, our coordinates
$w_1,...,w_{n-2},z$ give linear coordinates
on $\cal X$.  In these coordinates, the
function $A$ has the form
\begin{equation}
A(w_1,...,w_{n-2},z)=A'(w_1,...,w_{n-2})-cz\overline z.
\end{equation}
Here $c$ is some constant which depends on the cone
deficit of the added cone, and $A'$ is the area
form on the moduli space determined by the
list $\theta_1+\theta_2,\theta_3,...,\theta_n$.
Hence, if $A'$ has type $(1,n-3)$ then
$A$ has type $(1,n-2)$.

We have done the induction step but not the base
case.  We can do a reduction above unless $n=3$
or $n=4$ and all the cone points have the same
deficit.  The case $n=3$ is trivial -- the moduli
space is a single point.  When  $n=4$ and all
cone deficits are equal, the moduli space is
a finite cover of the modular surface,
which is modeled on $\C\H^1$.  So, the base 
cases work out.

\section{Some Definitions}
\label{defs}

\noindent
{\bf Stratified Manifolds:\/}
Let $\sqcup$ denote disjoint union.
A {\it complex hyperbolic stratified manifold\/} is a 
complete metric space
$X=X_0 \sqcup X_1 \sqcup X_2 \sqcup ...$
such that
\begin{itemize}
\item $X_0$ is connected and has finite volume.
\item $X_k$ is locally isometric to $\C\H^{n-k}$ for $k=0,1,2,...$
\item $X_{k+1} \subset {\rm closure\/}(X_k)$ for $k=0,1,2,...$.
\end{itemize}

\noindent
{\bf Fibered Cone Structure:\/}
We say that $X$ is a {\it fibered cone manifold\/}
if, for each $k$ and each $r \in X_k$, we have some
neighborhood $N_r$ of $r$ in $X_0$ with the
following structure. First
\begin{equation}
\label{fiber}
N_r=\bigcup_{s \in \Delta} F_s.
\end{equation}
Here $\Delta$ is an open metric disk about $r$ in $X_k$.
This is supposed to be a smooth fibration:
There is a fiber preserving diffeomorphism
$h: \Delta \times F_r \to N_r$.
We call $\Delta$ a {\it basic disk\/}.

Second, $\partial F_r$ is a spherical manifold,
and $F_r$ is foliated by geodesic arcs connecting
points on $\partial F_r$ to $r$.  These arcs
are all perpendicular to $\partial F_r$ at
their endpoints and they all have the same
length -- exactly the complex hyperbolic 
radius of the sphere
on which $\partial F_r$ is modeled.
\newline
\newline
{\bf Codimension 2 Conditions:\/}
$\C\H^n$ contains $\C\H^{n-1}$ as a totally geodesic
submanifold.  Let $Y$ denote the universal cover
of $\C\H^n-\C\H^{n-1}$.  
Let $\overline Y$ denote the metric completion of
$Y$.  The deck group extends to an action on
$\overline Y$. The projection $\pi: \overline Y \to \C\H^n$
is an infinite cyclic branched cover, branched over
$\C\H^{n-1}$.

For any $t>0$ there is an isometry
$I_t: \overline Y \to \overline Y$ which
rotates $\overline Y$ by an angle of $2 \pi t$ around
$\overline Y-Y$. A {\it simple cone manifold\/}
is a quotient of the form $\overline Y/I_t$ for
$t \in (0,1)$.  When $1/t \in \Z$, we call
$Y/I_t$ a {\it simple orbifold\/}.

We say that a complex hyperbolic stratified manifold has
codimension $2$ simple cone (respectively orbifold)
 singularities if every point
$p \in X_1$ has a neighborhood which is isometric to
a ball in a simple cone  manifold (respectively orbifold).
The isometry needs to take $X_1$ into the singular set,
and the dimensions are supposed to match up.

\section{Structure of the Completion}
\label{struct}

Now I'll revisit Thurston's paper and explain why
the completion of $\cal M$ has all the advertised
properties.  For ease of exposition, I'll work
with the labeled space $\cal M$, under the assumption
that the second angle condition simply does not
occur.  If the second angle condition does occur
then we first do the analysis in the labeled case
and then observe that the codimension $2$ orbifold
conditions emerge when we pass to the unlabeled
quotient space.  The point is that the relevant
cone angle gets cut in half.
\newline
\newline
{\bf Stratified Structure:\/}
Say that a {\it multi-list\/} is a subset
$\{\theta_{ij}\}$ of our list of deficits,
with $i \in \{1,...,k\}$,
such that
\begin{equation}
\alpha_i=\sum_j \theta_{ij}<2 \pi, \hskip 30 pt i=1,...,k
\end{equation}
We can take a sequence of points in $\cal M$ corresponding
to flat cone structures in which all the points 
corresponding to $\{\theta_{ij}\}$ for fixed $i$ coalesce.
The limit of this sequence is contained in a
point in the completion of $\cal M$ corresponding to
a stratum of codimension $\ell$.  Here
$\ell+3$ counts the number of cone points of
the limiting flat cone spheres.
This particular stratum is isometric to some
lower dimensional moduli space $\cal N$.
Moreover, any degeneration of structures in
$\cal M$ arises this way.  This gives the stratified
structure.
\newline
\newline
{\bf Fibered Cone Structure:\/}
Associated to $\cal N$ is the multi-list above.
Say that a cone point is {\it involved\/} if
it is one of the points corresponding to our
multi-list, and otherwise {\it uninvolved\/}.
Suppose we hold the uninvolved points fixed
and then coalesce the involved points.
This produces a point $r \in \cal N$.
There are $k$ Euclidean cones $C_1,...,C_k$
so that the $j$th cluster of involved points
coalseces down to the apex of $C_j$.
The point $r$ corresponds to the flat cone
structure defined by the apices of
the $C_j$ and the uninvolved points.
Fixing the uninvolved points and varying
the involved points gives $F_r$.

Why is $F_r$ totally geodesic?  As we did in
the previous section, we can choose local
linear coordinates so that the variables
$w_1,w_2,...$ describe the positions of
the uninvolved points and the positions
of the apices of the auxilliary cones,
and then variables $z_1,z_2,...$
describe the positions of the involved
points.  The points in $\C^{m-2}$
corresponding to $F_r$ comprise a complex
linear subspace.  So, when we projectivize,
we get some intersection of
$\C\H^{m-3}$ with a lower dimensional
complex projective space.  This gives
us a totally geodesic copy of a lower
dimensional complex hyperbolic place.

There is a natural foliation of $F_r$ into
arcs of geodesics.  If we start with one
flat cone structure corresponding to a point
in $F_r$, we can move the $j$th cluster of
involved points closer to the apex of $C_j$
by a homothety (i.e. similarity with
no twisting).  In terms of the coordinates
we just mentioned, we are simply replacing
$z_1,z_2,...$ with $rz_1,rz_2,...$ for
some real $r<1$.
We fix some $A-\epsilon$, where $A$ is the
area of the flat cone surface associated to
$r$ (before rescaling) and $\epsilon$ is some
small number.  If we restrict $F_r$ to 
points corresponding to structures having
area in $(A-\epsilon,A)$ then
$\partial F_r$ is a spherical manifold.

There is a natural diffeomorphism from $F_r$ to
a nearby fiber $F_s$ which comes from keeping
the involved points fixed relative to
the apices of the cones $C_1,...,C_k$
 and perturbing these apices and the
uninvolved points.  This map varies
smoothly with $s$ and gives rise to the
smooth fibration structure.
\newline
\newline
{\bf Codimension 2 Conditions:\/}
Now we consider the codimension $2$ strata.
The simplest case occurs when $k=1$ and $\{\theta_{ij}\}$
consists of just $2$ deficits whose sum is less than $2\pi$.
Each choice leads to a connected (real) codimension $2$ stratum.
A local analysis, as done in class, shows that a
neighborhood of each point on one of these strata
is isometric to the simple cone manifolds 
$Y/I_t$ discussed above.  The value of $t$ is
$2\pi-\theta_{11}-\theta_{12}$.
So, when Condition 1 on the deficits is satisfied, the corresponding
stratum is a codimension $2$ orbifold singularity.

A somewhat more subtle case occurs when
$\theta_{11}=\theta_{12}=\theta$.  The
analysis applied to the labeled moduli space
gives the angle around the stratum as
$2\pi-2\theta$.  However, as mentioned above,
when we pass to
the unlabeled moduli space, we are taking
a finite quotient which, in particular,
cuts this cone angle in half.
This case corresponds to Condition 2 on the deficits.
\newline
\newline
{\bf Finite Volume:\/} (This part seems to be done
just fine in Thurston's paper, and originally I
hadn't said anything about it in these notes.
But here I am adding some explanation.) 
Why does the space have finite volume?
The non-compact ends of the the space correspond to
paritions of the cone angles into two halves, each
of which sum to $\pi$.  There are finitely many
partitions, and you want to see that each partition leads
to a non-compact end with finite volume.  

Fix one of these partitions.
Outside a compact set,
the corresponding cone manifolds is a
``cigar'' -- a cylinder which has been 
capped off on either end.
If you fix, say, the minimum distance between the
cone points on the one end of the cigar and the
cone points on the other, then the set of all cone
structures realizing this minimum distance
is compact and hence has finite volume.

These ``fixed-minimum-distance'' 
 sets give a fibration of the
non-compact end. There is a natural operation
of inserting a cylinder in the middle (and rescaling the
area).  This insertion moves you from one fiber to another
one further out.  A local calculation shows that the insertion
of a cylinder of length $r$ decreases the volume of
the fiber by $C\exp(-r)$, for some constant $C$.
Hence, when you integrate over the fibers you get
finite volume. 

This argument is similar to the usual proof that a
cusped complex hyperbolic manifold has finite volume.

\section{Lattice Quotients}
\label{latq}

To finish the proof of Theorem \ref{lattice},
I'll prove the following result.

\begin{theorem}
\label{orbi}
\label{strat}
Let $X$ be a complex hyperbolic stratified manifold with
a fibered cone structure and
codimension $2$ orbifold conditions.  Then $X$ is a
lattice quotient. 
\end{theorem}

The rest of these notes are devoted to proving
Theorem \ref{strat}.
We make some basic definitions.

\begin{itemize}
\item Let $\widetilde X_0$ denote the universal
cover of $X_0$. 
\item  We have the developing map
$D: \widetilde X_0 \to \C\H^n$.
\item We have the holonomy homomorphism
$h: G=\pi_1(X_0) \to {\rm Isom\/}(\C\H^n)$.
\item Let $N \subset G$ be the kernel of the holonomy homomorphism.
\item  Let
$\widehat X_0=\widetilde X_0/N$.
\item $\widehat G=G/N$.
\item
The developing map factors through
a map $\widehat D: \widehat X_0 \to \C\H^n$.
\end{itemize}

Let $\widehat X$ denote the metric completion
of $\widehat X_0$. 

\begin{lemma} 
The developing map $\widehat D$
extends to $\widehat X$ and is a
distance non-increasing map.
Also, $\pi: \widehat X_0 \to X_0$
extends to a map $\pi: \widehat X \to X$.
\end{lemma}

\startproof
Choose a point $p \in \widehat X$.
There is a sequence $\{p_n\} \in \widehat X_0$
converging to $p$.  Define
$\widehat D(p)=\lim \widehat D(p_n)  \in \C\H^n$. 
This makes sense because $\widehat D$ is a local
isometry on $\widehat X_0$ and
hence distance non-increasing.  In particular
$\{\widehat D(p_n)\}$ is a Cauchy sequence.
If $\{q_n\}$ converges to $p$ as well then
$d(p_n,q_n) \to 0$.  But then
we must have 
$d(\widehat D(p_n),\widehat D(q_n)) \to 0$ as well.
Hence, $\widehat D(p)$ is well defined.
Since $\widehat D$ is distance non-increasing on
a dense subset of $\widehat X$, it is also
distance non-increasing on $\widehat X$.
The proof for $\pi$ is essentially the same.
\endproof

Define
\begin{equation}
\widehat X_k=\pi^{-1}(X_k), \hskip 30 pt k=1,2,3,...
\end{equation}

The next result is where we use the codimension $2$
orbifold conditions.

\begin{lemma}[Removable Singularities]
\label{good}
Every point $p \in \widehat X_1$ has a neighborhood
which is locally isometric to a ball in $\C\H^n$ and
the map $\widehat D$ gives such a local isometry.
\end{lemma}

\startproof
We first consider the picture in the space
$\widetilde X$. 
All the constructions above made for
$\widehat X$ also work for $\widetilde X$.
Let $C$ be a component of
$X_1$ and let $\widetilde C$ be a corresponding
component of $\widetilde X_1$.   The codimension $2$
cone manifold conditions tell us that
the map $\pi: \widetilde X \to X$ is an infinite
branched cover in a neighborhood of $\widetilde C$,
branched over $\widetilde C$.

We know that neighborhoods of points in $C$ are
locally isometric to balls in $Y/I_t$ for
some $t=1/k$.  Let $\beta$ be a small loop
in $\widetilde X_0$ which winds $k$ times
around $\widehat C$. The element (really conjugacy class
of elements) in the fundamental group $\pi_1(X_0)$
corresponding to $\beta$ has trivial holonomy, and
elements corresponding to loops winding fewer
times around have nontrivial holonomy.  For
this reason, $\widetilde X$ is isometric to a neighborhood of
$\overline Y/I_1=\C\H^n$ around $\widehat C$.
\endproof

The next lemma is where we use the fibered cone
manifold conditions.  This lemma seems obvious
at first glance, because 
$\dim(X_k)=\dim(X)-2k$.  However, this seems like
a slippery business.  So, I am going to spell out the proof
in a lot of detail.

\begin{lemma}[Dimension]
$\widehat D(\widehat X_k)$ has codimension at least $4$
for $k \geq 2$.
\end{lemma}

\startproof
Since $X_k$ has a compact exhaustion, $X_k$ is
covered by countably many basic disks. Hence,
it suffices to prove our result for
$\widehat \Delta=\pi^{-1}(\Delta)$, where
$\Delta$ is a basic disk.  Let $r$ be the
center of $\Delta$ and let $N=N_r$ be the
associated fibered neighborhood.

We claim that every point of
$\widehat p \in \widehat \Delta$ is an accumulation point
of a path component of $\pi^{-1}(N\cap X_0)$.
To see this, let $\{\widehat q_n\}$ be a Cauchy sequence
in $\widehat X_0$ converging to $\widehat p$.
The metric on $\widehat X_0$ is the path metric,
so we can find a path $\widehat \gamma_{mn}$ joining
$\widehat q_m$ to $\widehat q_n$ which is within
a factor of $2$ of the actual distance in
$\widehat X_0$ between $\widehat q_m$ and $\widehat q_n$.
Since $\pi$ does not increase distances,
$\{q_n\}$ is a Cauchy sequence in $X_0$ converging to $p$.
But then $q_n \in N$ for large $n$. Moreover
$\gamma_{mn} \subset N$ for large $m,n$.
But then the entire path $\widehat \gamma_{mn}$
stays in the same path component of $\pi^{-1}(N \cap X_0)$.
Hence, the tail end of our Cauchy sequence 
$\{\widehat q_n\}$ stays in the same path component.
This establishes the claim.

The space $\widehat X_0$ contains a countable dense set, and each
path component of
$\pi^{-1}(N \cap X_0)$ is an open set containing one point
in this dense set that is not contained in any of the others.
Therefore, there are only countable many components of
$\pi^{-1}(N \cap X_0)$.  In light of our claim above,
it suffices to prove, for an arbitrary
path component $\widehat A$ of
$\pi^{-1}(N \cap X_0)$, that 
$\widehat D(\widehat \Delta \cap {\rm closure\/}(\widehat A))$
has dimension $2n-2k$ in $\C\H^n$.  

Since $N$ is foliated by sets of the form $F_s$ for $s \in \Delta$,
we have the decomposition
\begin{equation}
\widehat A=\bigcup_{s \in \Delta} \widehat A_s,
\hskip 30 pt \widehat A_s=\pi^{-1}(F_s \cap X_0).
\end{equation}

Since $\pi$ is a local isometry on $\widehat X_0$, we
see that $\widehat A_s$ is foliated by geodesic arcs,
all of the same length, which meet
$\partial \widehat A_s$ at right
angles.  Here $\widehat A_s$ is a manifold
modeled on a complex hyperbolic sphere.
The foliating arcs give a retraction of $\widehat A$ onto
its $\partial A$. (Remember the the inner endpoints of
the foliating arcs are not part of $\widehat A$.)
Hence $\partial \widehat A$ is connected.
At the same time, the product structure on
$N \cap X_0$ gives a continuous retraction
from $N \cap X_0$ to each fiber
$F_s \cap X_0$. This continuous rectraction
lifts to a continuous rectraction from
$\partial \widehat A$ to
$\partial \widehat A_s$. Hence
$\partial \widehat A_s$ is connected.
Therefore, the image
$\widehat D(\partial \widehat A_s)$ is contained 
in a geodesic sphere $S_s$ in $\C\H^n$. Moreover,
$\widehat D$ maps the geodesic arcs foliating
$\widehat A_s$ to geodesic arcs perpendicular
to $S_s$ and pointing inward.  These
geodesic arcs all have the same length, so
they all meet at the center $c_s$ of $S_s$. 

Suppose that we have a Cauchy sequence
$\{\widehat q_n\}$ converging to some point
$\widehat s \in \widehat \Delta \cap {\rm closure\/}(\widehat A)$.
Then $\widehat q_n$ lies on some foliating
arc of some $\widehat A_{s_n}$.
Since there is a minimum positive distance
between $s$ and any fiber $F_t$ with $t \not = s$,
we must have $s_n \to s$.  
Moreover, when $n$
is large, $\widehat q_n$ lies almost all the
way at the inner end of the foliating arc.
Hence, the distance from
$\widehat D(\widehat q_n)$ to $c_s$ tends to $0$
as $n$ tends to $\infty$. 
Hence
\begin{equation}
\widehat D(\widehat \Delta \cap {\rm closure\/}(\widehat A)) \subset
Y=\bigcup_{s \in \Delta} c_s
\end{equation}
Given the smooth nature of the fibration, the point $c_s$
varies smoothly with $s \in \Delta$.
This shows that $Y$ is a
smooth manifold of dimension $2n-2k$.
\endproof

\begin{lemma}
\label{remove}
$\widehat D(\widehat X)=\C\H^n$.
\end{lemma}

\startproof
From Lemma \ref{good}, the map $\widehat D$ is a
local isometry from $\widehat X_0 \cup \widehat X_1$ to 
$\C\H^n$.  Suppose $\widehat D: \widehat X \to \C\H^n$ is
not onto.
Let $q \in \C\H^n-\widehat D(\widehat X)$.
Pick $p \in \widehat D(\widehat X_0 \cup \widehat X_1)$ and consider the
geodesic $\gamma$ connecting $p$ to $q$.  Choosing $p$
generically and using the Dimension Lemma,
we can arrange that $\gamma$ does not intersect
$\widehat D(\widehat X_k)$ for $k \geq 2$.

Let $\widehat p$ be some pre-image of $p$ in $\widehat X_0 \cup \widehat X_1$.
There is some initial geodesic segment $\widehat \alpha$ which
$\widehat D$ carries to the initial portion of $\gamma$
emanating from $p$.  The geodesic $\widehat \gamma$ extending
$\widehat \alpha$ lies entirely in $\widehat X_0 \cup \widehat X_1$, by
construction.  But then $\widehat D$ is defined on
all of $\widehat \gamma$ and in particular
$q \in \widehat D(\widehat \gamma) 
\subset \widehat D(\widehat X_0 \cup \widehat X_1)$.
This is a contradiction.
\endproof

\begin{lemma}
$\widehat D$ is injective on $\widehat X_0 \cup \widehat X_1$.
\end{lemma}

\startproof
By Lemma \ref{remove},
the map $\widehat D: \widehat X_0 \cup \widehat X_1 \to \C\H^n$
is a local isometry and therefore a covering map
of its image.
But, by the Dimension Lemma and Lemma \ref{remove}, the image
$\widehat D(\widehat X_0 \cup \widehat X_1)$ is
everything but a set of codimension at least $4$. Hence
$\widehat D(\widehat X_0 \cup \widehat X_1)$ is
simply connected.  But then our covering map
must be injective.
\endproof

$\widehat D$ is a global isometry from
$\widehat X_0 \cup \widehat X_1$ to an open dense subset
of $\C\H^n$.  So, we can identify
$\widehat X_0 \cup \widehat X_1$ with an open dense
subset of $\C\H^n$. We make this identification.
Let $\Gamma=\widehat G$.
Under our identification, $\C\H^n$ is the
metric completion of $\widehat X_0 \cup \widehat X_1$.
Hence $\widehat X=\C\H^n$.  But then $\Gamma$
acts isometrically on $\C\H^n$.
The action is discrete and co-finite
because $\widehat X_0/\Gamma=X_0$ has finite volume
and nonempty interior.  
Hence $\Gamma$ is a lattice.  Since
$X_0$ is dense in $\C\H^n/\Gamma$, and
$\C\H^n/\Gamma$ has finite volume, we see that
$\C\H^n/\Gamma$ is the metric completion of
$X_0$.  Hence $X=\C\H^n/\Gamma$.  This completes the proof.

\section{Proof of Theorem \ref{tri1}}
\label{triproof1}

I'll prove Theorem \ref{tri1}
through a series of smaller results.

\begin{lemma}
If ${\cal M\/}$ is special, the set of
triangulation points is dense in $\cal M$.
\end{lemma}

\startproof
Let's look at the local coordinates we get when we have
a triangulation.  We take some embedded spanning tree and
associate the coordinates as above.  When we cut along
the spanning tree and look at the resulting pseudo-polygon,
we can develop it into $\C$ so that the vertices lie in
${\bf Eis\/}$.  Moreover, the unit complex numbers
$\{u_i\}$ relating pairs of coordinates (on edges which
get glued together) also belong to ${\bf Eis\/}$.  In short,
all the coordinates lie in ${\bf Eis\/}$.  Conversely,
if we choose sufficiently nearby
coordinates in ${\bf Eis\/}$, we get a
triangulation.

Now, if we have any flat cone sphere corresponding to a
point in $\bf M$, we can scale it up so that it has
enormous coordinates with respect to some spanning
tree, and then we can find nearby coordinates in
${\bf Eis\/}$ which just differ by at most $2$ units
from the original coordinates.  When we scale back
down to (say) unit area, the original structure
and the nearby triangulation point are extremely
close.  Hence the triangulation points are dense
in $\bf M$.
\endproof

\begin{lemma}
If ${\cal M\/}$ is special, then the 
completion of $\cal M$ is a
lattice quotient.
\end{lemma}

\startproof
We just have to verify the deficit conditions.
We have $\theta_i=k_i \pi/3$ for some $k_i \in \{1,2,3\}$.
If $\theta_i+\theta_j<2\pi$ and $\theta_i \not = \theta_j$ then
$2\pi-\theta_i-\theta_j=k_{ij}\pi/3$
for some $k_{ij}=1,2,3$.  Hence, the first condition on
the deficits holds.
If $\theta_i=\theta_j$ and $\theta_i<\pi$ then,
again $2\pi$ is an integer multiple of $\pi-\theta_i$.
Now we apply Theorem \ref{lattice}.
\endproof

Let $\Gamma$ be the lattice
such that the completion of ${\cal M\/}$ is
$\C\H^{m-3}/\Gamma$.  Now,
$\Gamma$ acts as a group of matrices on
$\C^{m-2}$ and preserves some Hermitian
form $A$ of type $(1,m-2)$.  To get actual
matrices, we need to choose some linear
coordinates on $\C^{m-2}$.  We choose
the coordinates coming from the embedded
spanning trees.
This gives 
$2m-2$ variables, but we choose $m-2$ independent ones.

\begin{lemma}
With respect to the coordinates coming from the 
embedded spanning trees, the
entries of elements of $\Gamma$ all lie in
${\bf Eis\/}$.
\end{lemma}

\startproof
Suppose we start with a closed loop in $\cal M$.
We develop $\cal M$ into $\C\H^{m-2}$ along
this loop and then take the holonomy.  This
gives us some element of $\Gamma$, and all
elements of $\Gamma$ arise this way.

We can break our loop into finitely many
segments, such that each segment is contained
in a single spanning tree coordinate chart on $\cal M$.  
As we move from segment to segment, we make
some linear change of coordinates.  The
element of $\Gamma$ is the product of these coordinate-change
matrices.

Now, when we compute the coordinate change matrices, we
can compute them with respect to triangulation
points, because the triangulation points are dense.
But, from the description of the coordinate
changes given in \S \ref{loc} we see that each
coordinate on the new spanning tree is a complex
linear combination of the old coordinates, where
the coefficients of the linear combination
lie in ${\bf Eis\/}$.  Hence, the coordinate
change matrices have entries in ${\bf Eis\/}$.
Hence, so does the product of these matrices.
\endproof

\begin{lemma}
$\Gamma$ preserves a Hermitian form $H$ of
type $(1,9)$ which is defined over ${\bf Eis\/}$.
\end{lemma}

\startproof
Choose some point $p \in \cal M$ and
let $\Sigma$ be the corresponding 
flat cone sphere.
Let $\tau$ be some embedded spanning tree
on $\Sigma$.  When we develop out $\Sigma-\tau$
into the plane, we get $2n-2$ coordinates
which, in pairs, are related by complex numbers
$u_1,...,u_{n-1}$.  But, due to the values of
the cone angles, these numbers are all
$6$th roots of unity.  They all belong to
${\bf Eis\/}$.  So, when we triangulate
the pseudo-polygon $\Sigma-\tau$, the other
labels are complex linear combinations of
the original variables, with coefficients
in ${\bf Eis\/}$.
\endproof

\noindent
{\bf Proof of Theorem \ref{tri1}:\/}
The lattice $\Gamma$ corresponds to the moduli
space which contains the regular icosahedral
tiling.  This lattice acts on
$\C\H^9$ because there are $12$ deficits
on the list and $9=12-3$.
If we coalesce various of the cone
points corresponding to the regular icosahedron,
we can achieve a deficit list corresponding
to every other type of triangulation. (For
instance, if we coalesce the points in pairs
we get triangulations having the same deficit
list as the list produced by the octahedron.)
But this means that every special moduli
space is some stratum of $\C\H^9/\Gamma$.
From what we have already seen, 
$\Gamma$ is defined over ${\bf Eis\/}$.
\endproof

\section{Proof of Theorem \ref{tri2}}
\label{triproof2}

Let $\Gamma$ be the lattice from
Theorem \ref{tri1}. 
We fix some point in $\cal M$, say
the structure corresponding to the
regular icosahedron.  We also fix
some embedded spanning tree $\tau$
relative this structure.
There is some open cone $\cal C$ in
$\C^{1,9}$ such that points in
$\cal C$ correspond, via coordinates
on $\tau$, to some open set in
$\cal M$.   We define the Hermitian
form $H=4 \sqrt 3 A$ with respect to $\tau$.
The elements of $\Gamma$ preserve
both $H$ and ${\bf Eis\/}^{1,9}$, even
though they typically move the
cone $\cal C$ off itself.

\begin{lemma}
$H$ is defined by a metrix with entries in
${\bf Eis\/}$.
Given a positive
vector $V \in {\cal C\/} \cap {\bf Eis\/}^{1,9}$,
the norm $H(V,V)$ computes the $4$ times
the area of the flat cone sphere.
\end{lemma}

\startproof
When we work out the formula for $H$ with
respect to $\tau$, we see that it
just involves expressions of the form
$\sqrt 3 i(z \overline w - w\overline z)$, where
$z$ and $w$ are complex linear combinations
of the coordinates with coefficients
in ${\bf Eis\/}$.  This easily implies
that $H$ is defined over ${\bf Eis\/}$.

For the second statement, we observe that
a unit equilateral triangle has area
$\sqrt 3/4$. Hence the quantity $H(V,V)$,
which records $4 \sqrt 3$ times the area, counts
$3$ times the number of triangles.
\endproof

\noindent
{\bf From Triangulations to Vectors:\/}
For any pair
$(\Sigma,\sigma)$, where $\Sigma$ is a flat
cone sphere and $\sigma$ is an embedded
spanning tree, there is some finite
sequence of coordinate changes whose
composition allows us to express
the $\sigma$-coordinates as
$\tau$-coordinates.
What we have is a finite sequence
$(\Sigma_j,\sigma_j)$  and
a finite sequence ${\cal C\/}_j$
of cones such that each point in
${\cal C\/}_j$ corresponds, via coordinates
on $\sigma_j$ to a flat cone sphere.
Here $j=0,...,k$ and $\sigma_0=\sigma$ and
$\sigma_k=\tau$.
The cones ${\cal C\/}_j$ and
${\cal C\/}_{j+1}$ overlap, and there is
some matrix $M_{j+1}$, defined over ${\bf Eis\/}$, which
expresses the coordinate changes on the
overlap.   
The product of the matrices $M=M_k...M_1$
expresses the $\sigma$-coordinates
in terms of the $\tau$ coordinates,
even though the $\tau$ coordinates may
not lie in the cone ${\cal C\/}={\cal C\/}_k$.
This does not bother us. The important point is that
these coordinate changes preserve both ${\bf Eis\/}^{1,9}$
and $H$.

Now, suppose we were to take a different
sequence $(\Sigma_j',\sigma_j')$ for $j=1,...,\ell$,
with $\sigma_0'=\sigma$ and
$\sigma_{\ell}'=\tau$.  This would give us
sequence of cones and matrices, and hence
a new coordinate change.  In this case,
we would apply the matrix
$M'=M_{\ell}'...M_1'$ to the $\sigma$-coordinates.
The matrix $M' \circ M^{-1}$ is the result of
doing a ``loop of coordinate changes'' starting
and ending at $(\Sigma,\tau)$.  Hence, this matrix
belongs to $\Gamma$.  In short, any two of our
coordinate changes differ by the action of
$\Gamma$.  In other words, given $\sigma$-coordinates,
there is a canonical point in $\C^{9,1}/\Gamma$ that
represents the $\tau$-coordinates mod $\Gamma$.

So, if we start with a triangulation of the
sphere, we get a point in ${\bf Eis\/}^{1,9}$ 
relative to some embedded spanning tree.
We then make a coordinate change and get
a well-defined positive vector in 
${\bf Eis\/}^{1,9}/\Gamma$.
\newline
\newline
{\bf From Vectors to Triangulations:\/}
Conversely, if we have a positive vector
in ${\bf Eis\/}^{1,9}/\Gamma$ we take
some representative vector $V \in {\bf Eis\/}^{1,9}$
and then interpret $V$ as giving coordinates relative
to our preferred tree $\tau$, even though
these coordinates might not lie in the cone
$\cal C$.  We then make
a finite string of coordinate changes
until we arrive at a new vector 
$W \in {\bf Eis\/}^{1,9}$ giving coordinates
relative to a spanning tree $\sigma$ which
is embedded on the flat cone sphere $\Sigma$ corresponding
to $[W]$.  This gives us a triangulation
of the sphere.  

This triangulation is
independent of our choice of coordinate
change, and also independent of the
choice of $V$.  If we make the
construction twice, the two triangulations
on $\Sigma$ have $\sigma$-coordinates
which differ by a loop of coordinate
changes, as above, starting from and ending at
$(\Sigma,\sigma)$.  These coordinate
changes do nothing to the triangulation.
\newline

The two halves of our construction are
inverses of each other, so we get a bijection
between the advertised sets.
Since $\Gamma$ is constructed out of the
kind of sequences of coordinate changes
discussed above, $\Gamma$ preserves the
Hermitian form $H$.  As we have already
mentioned, $H(V,V)$ counts
$3$ times the number of triangles.

\section{References}

Here are some additional refrences:
\begin{itemize}
\item Curt McMullen's paper
{\it The Gauss-Bonnet Theorem for Cone Manifolds and Volumes of Moduli Spaces\/}.
works out a general theory of cone manifolds which adds details
to Thurston's description, especially a
prime factorization theorem for cone manifolds.

\item M. Weber's thesis 1993 Bonn thesis
{\it Fundamentalbereiche komplex hyperbolischer Fl \"achen.\/}
works through Thurston's construction using
star-shaped spanning trees.  I don't know how well this matches
what I do above.

\item John Parker's paper
{\it J.R. Parker, Cone metrics on the sphere and Livné's lattices\/}
(Acta Mathematica 196 (2006) 1-64)
Works through Thurston's construction explicitly for Livn\'e's lattices,
including building fundamental polyhedra using the Poincar\'e theorem.
The Livn\'e's lattices are special cases, corresponding to moduli
spaces with $5$ cone points.

\item An upcoming paper by 
John Parker and Richard K Boadi,
{\it Mostow's lattices and cone metrics on the sphere\/}
(Advances in Geometry)
does the same thing as Parker's earlier
paper but for some of Mostow's lattices.

\end{itemize}

\end{document}